\documentclass[a4paper,11pt]{amsart}
\usepackage{graphicx}
\usepackage{mathptmx}
\usepackage{amsmath}
\usepackage{amssymb}
\usepackage{enumitem}
\usepackage{xcolor}

\newmuskip\pFqmuskip

\newcommand*\pFq[6][8]{%
  \begingroup % only local assignments
  \pFqmuskip=#1mu\relax
  \mathcode`=\string"8000
  \begingroup\lccode`\~=`\,
  \lowercase{\endgroup\let~}\pFqcomma
  F^{#2}_{#3}{\left(\genfrac..{0pt}{}{#4}{#5}\bigg|#6\right)}%
  \endgroup
}
\newcommand{\pFqcomma}{\mskip\pFqmuskip}

\newtheorem{theorem}{Theorem}

\newtheorem{corollary}[theorem]{Corollary}

\newtheorem{remark}[theorem]{Remark}

\begin{document}
\title[]{Identities involving degenerate harmonic and degenerate hyperharmonic numbers}
\author{Hye Kyung Kim$^{1}$}
\address{Department Of Mathematics Education, Daegu Catholic University, Gyeongsan 38430, Republic of Korea}
\email{hkkim@cu.ac.kr}

\author{Dae San Kim$^{2}$}
\address{Department of Mathematics, Sogang University, Seoul 121-742, Republic of Korea}
\email{dskim@sogang.ac.kr}

\author{Taekyun Kim$^{3,*}$}
\address{Department of Mathematics, Kwangwoon University, Seoul 139-701, Republic of Korea}
\email{tkkim@kw.ac.kr}

\thanks{* is corresponding author}

\begin{abstract}
Harmonic numbers have been studied since antiquity, while hyperharmonic numbers were intoduced by Conway and Guy in 1996. The degenerate harmonic numbers and degenerate hyperharmonic numbers are their respective degenerate versions.
The aim of this paper is to further investigate some properties, recurrence relations and identities involving the degenerate harmonic and degenerate hyperharmonic numbers in connection with degenerate Stirling numbers of the first kind, degenerate Daehee numbers and degenerate derangements.
\end{abstract}

 \maketitle

\markboth{\centerline{\scriptsize Identities involving degenerate harmonic and degenerate hyperharmonic numbers}}
{\centerline{\scriptsize H. K. Kim, D. S. Kim, and T. Kim}}

\section{Introduction}

\medskip

%%%%%%%%%%%%%%%%%%%%%%%%%%%%%%%%%%%%%%%%%%%%%%%%%%%%%%%%%%%%%%%%%%%%%%%%%%%%%%%%%%%%%0104
%1 

In recent years, various degenerate versions of many special numbers and polynomials have beem studied and yielded a lot of fascinating and fruitful results (\text{see \cite{5, 6, 7, 8, 9, 10, 11, 12}} and the references therein), which began with Carlitz's work on the degenerate Bernoulli and degenerate Euler numbers (see [2]). It is worthwhile to mention that these explorations for degenerate versions are not limited to polynomials and numbers but also extended to transcendental functions, like gamma functions (\text{see \cite{9,10}}). It is also remarkable that the $\lambda$-umbral calculus and $\lambda$-$q$-umbral calculus were introduced as degenerate versions of the umbral calculus and the $q$-umbral calculus, respectively (\text{see \cite{6, 11}}). As it turns out, the $\lambda$-umbral calculus and $\lambda$-$q$-umbral calculus are more convenient than the umbral calculus and the $q$-umbral calculus when dealing with degenerate Sheffer polynomials and degenerate $q$-Sheffer polynomials. \par
The aim of this paper is to further investigate some properties, recurrence relations and identities involving the degenerate harmonic numbers (see \eqref{eq06}) and the degenerate hyperharmonic numbers (see \eqref{eq07}, \eqref{eq08}) in connection with degenerate Stirling numbers of the first kind, degenerate Daehee numbers and degenerate derangements. The degenerate harmonic numbers and degenerate hyperharmonic numbers are respectively degenerate versions of the harmonic numbers and the hyperharmonic numbers, of which the latter are introduced in \text{\cite{4}}. \par
The outline of this paper is as follows. In Section 1, we recall the degenerate exponentials and the degenerate logarithms. We remind the reader of the harmonic numbers, and of the hyperharmonic numbers together with their explicit expression due to Conway and Guy (\text{see \cite{4}}). Then we recall their degenerate versions, namely the degenerate harmonic numbers, and the degenerate hyperharmonic numbers together with their explicit expression (\text{see \cite{7,8}}). We also mention the recently introduced degenerate Stirling numbers of the first kind and the degenerate Daehee numbers of order $r$.
Section 2 is the main result of this paper. We obtain an expression of the degenerate hyperharmonic numbers of order $r$ in terms of the same numbers of lower orders in Theorem 1. We express the Daehee numbers in terms of the degenerate harmonic numbers and of the degenerate hyperharmonic numbers, respectively in Theorem 2 and Theorem 3.
In Theorem 4, the degenerate harmonic numbers are represented in terms of the degenerate hyperharmonic numbers of order $r$. In Theorem 5, the degenerate Daehee numbers are represented in terms of the degenerate Daehee numbers of order $r-1$ and of the degenerate hyperharmonic numbers. We derive a simple relation between the degenerate hyperharmonic numbers and the degenerate Daehee numbers in Theorem 6. We deduce an identity involving the degenerate hyperharmonic numbers and the degenerate derangements in Theorem 7. The degenerate Daehee numbers are expressed in terms of the degenerate Stirling numbers of the first kind in Theorem 8. Finally, we get an identity involving the degenerate Stirling numbers of the first kind and the degenerate harmonic numbers in Theorem 9. \par

\vspace{0.1in}
For any nonzero $\lambda\in\mathbb{R}$, the degenerate exponential functions are defined by
\begin{equation}\label{eq01}
\begin{split}
e_\lambda^x(t)=(1+\lambda t)^{\frac{x}{\lambda}}=\sum_{n=0}^\infty(x)_{n,\lambda}\frac{t^n}{n!},\quad e_\lambda(t)=e_\lambda^1(t),\quad (\text{see \cite{2, 8}}),
\end{split}
\end{equation}
where
\begin{equation*}
\begin{split}
(x)_{0,\lambda}=1,\ \ (x)_{n,\lambda}=x(x-\lambda)\cdots (x-(n-1)\lambda), \ (n\geq 1),\quad (\text{see \cite{8}}).
\end{split}
\end{equation*}
Let $\log_\lambda t$ be the compositional inverse of $e_\lambda(t)$ with $e_\lambda(\log_\lambda t)=\log_\lambda e_\lambda(t)=t$. It is called the degenerate logarithm and is given by
\begin{equation}\label{eq02}
\begin{split}
\log_\lambda(1+t)=\sum_{k=1}^\infty \frac{\lambda^{k-1}(1)_{k,\frac{1}{\lambda}}}{k!}t^k=\frac{1}{\lambda}((1+t)^\lambda-1),\quad (\text{see \cite{5}}).
\end{split}
\end{equation} \par
The harmonic numbers are given by
\begin{equation}\label{eq03}
\begin{split}
H_0=0, \ \ H_n=1+\frac{1}{2}+\cdots+\frac{1}{n}, \quad (n\in\mathbb{N}),\quad \ (\text{see \cite{3, 4, 16}}).
\end{split}
\end{equation}
In 1996, Conway and Guy introduced the hyperharmonic numbers $H_n^{(r)}$ of order $r$,\,$(n,r\geq0)$, which are given by
\begin{equation}\label{eq04}
\begin{split}
H_0^{(r)}=0,\,\,(r \ge 0),\,\, H_n^{(0)}=\frac{1}{n},\,\,(n \ge 1),\,\, H_n^{(r)}=\sum_{k=1}^nH_k^{(r-1)}, \,\, (n, r\ge 1),\quad \ (\text{see \cite{4}}).
\end{split}
\end{equation}
Thus, by \eqref{eq04}, we get
\begin{equation}\label{eq05}
\begin{split}
H_n^{(r)}=\binom{n+r-1}{n}(H_{n+r-1}-H_{r-1}),\quad (r\geq1),\quad \ (\text{see \cite{4}}).
\end{split}
\end{equation}
Recently, the degenerate harmonic numbers are defined by
\begin{equation}\label{eq06}
\begin{split}
H_{0,\lambda}=0, \ H_{n,\lambda}=\sum_{k=1}^n\frac{1}{\lambda}\binom{\lambda}{k}(-1)^{k-1}, \quad (n\geq1),\quad \ (\text{see \cite{8}}).
\end{split}
\end{equation}
Note that $\lim_{\lambda\rightarrow 0}H_{n,\lambda}=H_n$.
The degenerate hyperharmonic numbers $H_{n,\lambda}^{(r)}$ of order $r$,\,$(n,r \ge 0)$,\, are defined by
\begin{equation}\label{eq07}
\begin{split}
H_{0,\lambda}^{(r)}=0,\,\,(r \ge 0),\,\, H_{n,\lambda}^{(0)}=\frac{1}{\lambda}\binom{\lambda}{n}(-1)^{n-1},\,\, (n \ge 1),\,\,H_{n,\lambda}^{(r)}=\sum_{k=1}^nH_{k,\lambda}^{(r-1)},\,\, (n,r\ge 1),\quad (\text{see \cite{7}}).
\end{split}
\end{equation}
We see from \eqref{eq06} and \eqref{eq07} that $H_{n,\lambda}^{(1)}=H_{n,\lambda}$.
From \eqref{eq07}, we note that
\begin{equation}\label{eq08}
\begin{split}
H_{n,\lambda}^{(r)}=\frac{(-1)^{r-1}}{\binom{\lambda-1}{r-1}}\binom{n+r-1}{n}(H_{n+r-1,\lambda}-H_{r-1,\lambda}),\quad (\text{see \cite{7}}),
\end{split}
\end{equation}
where $n,\ r$ are positive numbers. Here we observe from \eqref{eq05} and \eqref{eq08} that $\lim_{\lambda\rightarrow 0}H_{n,\lambda}^{(r)}=H_n^{(r)}$. \par
In \cite{5}, the degenerate Stirling numbers of the first kind are defined by
%3
\begin{equation}\label{eq09}
\begin{split}
(x)_n=\sum_{k=0}^nS_{1,\lambda}(n,k)(x)_{k,\lambda},\quad (n\geq0), \quad (\text{see \cite{5, 8}}),
\end{split}
\end{equation}
where $(x)_0=1,\ (x)_n=x(x-1)\cdots(x-n+1)$, $(n\geq1)$. \par
For $r\in\mathbb{N}$, the degenerate Daehee numbers of order $r$ are defined by 
\begin{equation}\label{eq10}
\begin{split}
\bigg(\frac{\log_\lambda(1+t)}{t}\bigg)^r=\sum_{n=0}^\infty D_{n,\lambda}^{(r)}\frac{t^n}{n!},\quad (\text{see \cite{11}}).
\end{split}
\end{equation}
In particular, for $r=1$, $D_{n,\lambda}=D_{n,\lambda}^{(1)}$ are called the degenerate Daehee numbers \par

\section{Identities involving degenerate harmonic and degenerate hyperharmonic numbers}
%4
From \eqref{eq06} and \eqref{eq07}, we note that
\begin{equation}\label{eq11}
\begin{split}
-\frac{\log_\lambda(1-t)}{(1-t)}=\sum_{n=1}^\infty H_{n,\lambda}t^n,\quad (\text{see \cite{7}}),
\end{split}
\end{equation}
and
\begin{equation}\label{eq12}
\begin{split}
-\frac{\log_\lambda(1-t)}{(1-t)^r}=\sum_{n=1}^\infty H_{n,\lambda}^{(r)}t^n,\quad (\text{see \cite{7}}),
\end{split}
\end{equation}
where $r$ is a nonnegative integer.

By \eqref{eq12}, we get
\begin{equation}\label{eq13}
\begin{split}
\sum_{n=1}^\infty H_{n,\lambda}^{(r-1)}t^n&=-\frac{\log_\lambda(1-t)}{(1-t)^r}(1-t)=\sum_{n=1}^\infty H_{n,\lambda}^{(r)}t^n(1-t)\\
&=\sum_{n=1}^\infty H_{n,\lambda}^{(r)}t^n-\sum_{n=1}^\infty H_{n,\lambda}^{(r)}t^{n+1}=\sum_{n=1}^\infty (H_{n,\lambda}^{(r)}-H_{n-1,\lambda}^{(r)})t^n.
\end{split}
\end{equation}
By comparing the coefficients on both sides of \eqref{eq13}, we get
\begin{equation}\label{eq14}
H_{n,\lambda}^{(r)}=H_{n-1,\lambda}^{(r)}+H_{n,\lambda}^{(r-1)}.
\end{equation}

For $1\leq s \leq r$, by \eqref{eq12}, we get
\begin{equation}\label{eq16}
\begin{split}
\sum_{n=1}^\infty H_{n,\lambda}^{(r)}t^n&=-\frac{\log_\lambda(1-t)}{(1-t)^r}=-\frac{\log_\lambda(1-t)}{(1-t)^{r-s}}\frac{1}{(1-t)^s}\\
&=\sum_{l=1}^\infty H_{l,\lambda}^{(r-s)}t^l \sum_{k=0}^\infty \binom{k+s-1}{k}t^k\\
&=\sum_{n=1}^\infty \sum_{l=1}^n H_{l,\lambda}^{(r-s)}\binom{n-l+s-1}{s-1}t^n.
\end{split}
\end{equation}
By comparing the coefficients on both sides of \eqref{eq16}, we get
\begin{equation}\label{eq17}
\begin{split}
H_{n,\lambda}^{(r)}=\sum_{l=1}^nH_{l,\lambda}^{(r-s)}\binom{n-l+s-1}{s-1},
\end{split}
\end{equation}
where $r,\ s\in \mathbb{Z}$ with $1\leq s \leq r$. 
In particular, for $r=s$, we have
\begin{equation}\label{eq18}
\begin{split}
H_{n,\lambda}^{(r)}=\sum_{l=1}^n H_{l,\lambda}^{(0)}\binom{n-l+r-1}{r-1}=\sum_{l=1}^n\frac{1}{\lambda}\binom{\lambda}{l}(-1)^{l-1}\binom{n-l+r-1}{r-1}.
\end{split}
\end{equation}
Therefore, by \eqref{eq17} and \eqref{eq18}, we obtain the following theorem.

\begin{theorem}
For $r,\ s \in \mathbb{Z}$ with $1\leq s \leq r$, we have
\begin{equation*}
\begin{split}
H_{n,\lambda}^{(r)}=\sum_{l=1}^nH_{l,\lambda}^{(r-s)}\binom{n-l+s-1}{s-1},
\end{split}
\end{equation*}
and
\begin{equation*}
\begin{split}
H_{n,\lambda}^{(r)}=\sum_{l=1}^n\frac{1}{\lambda}\binom{\lambda}{l}(-1)^{l-1}\binom{n-l+r-1}{r-1}.
\end{split}
\end{equation*}

\end{theorem}

From \eqref{eq11} and \eqref{eq14}, we note that
\begin{equation}\label{eq19}
\begin{split}
\sum_{n=0}^\infty D_{n,\lambda}\frac{t^n}{n!}&=\frac{\log_\lambda(1+t)}{t}=\frac{\log_\lambda(1+t)}{1+t}\frac{1+t}{t}\\
&=\bigg(\sum_{k=1}^\infty (-1)^{k+1}H_{k,\lambda}t^k\bigg)\bigg(1+\frac{1}{t}\bigg)\\
&=\sum_{n=1}^\infty(-1)^{n+1}H_{n,\lambda}t^n+\sum_{n=0}^\infty (-1)^n H_{n+1,\lambda}t^n\\
&=1+\sum_{n=1}^\infty(-1)^n(H_{n+1,\lambda}-H_{n,\lambda})t^n.
\end{split}
\end{equation}

%7
Therefore, by comparing the coefficients on both sides of \eqref{eq19}, we have the following theorem.

\begin{theorem}
For $n\geq0$, we have
\begin{equation*}
\begin{split}
D_{0,\lambda}=1,\,\, D_{n,\lambda}=(-1)^n n!(H_{n+1,\lambda}-H_{n,\lambda}),\,\,(n \ge 1).
\end{split}
\end{equation*}
\end{theorem}

From \eqref{eq12}, we note that
\begin{equation}\label{eq20}
\begin{split}
\sum_{n=0}^\infty D_{n,\lambda}\frac{t^n}{n!}&=\frac{\log_\lambda(1+t)}{t}=\frac{\log_\lambda(1+t)}{t (1+t)^r}(1+t)^r\\
&=\sum_{k=0}^\infty H_{k+1,\lambda}^{(r)}(-1)^kt^k\sum_{l=0}^\infty \binom{r}{l}t^l\\
&=\sum_{n=0}^\infty \bigg(\sum_{k=0}^nH_{k+1,\lambda}^{(r)}\binom{r}{n-k}(-1)^k\bigg)t^n.
\end{split}
\end{equation}

Therefore, by \eqref{eq20}, we obtain the following theorem

\begin{theorem}
For $n \ge 0$, we have
\begin{equation*}
\begin{split}
D_{n,\lambda}=n!\sum_{k=0}^nH_{k+1,\lambda}^{(r)}\binom{r}{n-k}(-1)^k.
\end{split}
\end{equation*}
\end{theorem}

%8
Now, we observe from \eqref{eq02} that
\begin{equation}\label{eq21}
\sum_{n=0}^\infty D_{n,\lambda}\frac{t^n}{n!}=\frac{\log_\lambda (1+t)}{t} \\
=\sum_{n=1}^\infty \binom{\lambda}{n} \frac{1}{\lambda}t^{n-1}=\sum_{n=0}^\infty \binom{\lambda}{n+1}\frac{1}{\lambda}t^n.
\end{equation}
Thus, by \eqref{eq21}, we get
\begin{equation}\label{eq22}
\begin{split}
D_{n,\lambda}=n!\frac{1}{\lambda}\binom{\lambda}{n+1}=\frac{(\lambda-1)_{n}}{n+1},\quad (n\geq0).
\end{split}
\end{equation} \par
From \eqref{eq11}, we have
\begin{equation}\label{eq23}
\begin{split}
\sum_{n=1}^\infty H_{n,\lambda}t^n&=-\frac{\log_\lambda (1-t)}{1-t}=-\frac{\log_\lambda (1-t)}{t}\frac{t}{1-t} \\
&=\sum_{l=0}^\infty D_{l,\lambda}(-1)^l \frac{t^l}{l!}\sum_{m=1}^\infty t^{m} \\
&=\sum_{n=1}^\infty \bigg(\sum_{l=0}^{n-1}D_{l,\lambda}\frac{(-1)^l}{l!}\bigg)t^n.
\end{split}
\end{equation}
Thus, by Theorem 3 and \eqref{eq23}, we get
\begin{equation}\label{eq24}
\begin{split}
H_{n,\lambda}&=\sum_{l=0}^{n-1} D_{l,\lambda} \frac{(-1)^l}{l!}=\sum_{l=0}^{n-1}\frac{(-1)^l}{l!}l! \sum_{k=0}^l H_{k+1,\lambda}^{(r)} \binom{r}{l-k}(-1)^k \\
&=\sum_{l=0}^{n-1}\sum_{k=0}^l (-1)^{k+l} H_{k+1,\lambda}^{(r)}\binom{r}{l-k}, \quad{(n\geq1)}.
\end{split}
\end{equation}
Therefore, by \eqref{eq24}, we obtain the following theorem.

\begin{theorem}
For $n\geq1$, we have
\begin{equation*}
\begin{split}
H_{n,\lambda}=\sum_{l=0}^{n-1}\sum_{k=0}^l(-1)^{k+l}\binom{r}{l-k} H_{k+1,\lambda}^{(r)}.
\end{split}
\end{equation*}
\end{theorem}

By \eqref{eq10}, we get
\begin{equation}\label{eq25}
\begin{split}
\sum_{n=0}^\infty D_{n,\lambda}^{(r)} \frac{t^n}{n!}&=\bigg(\frac{\log_\lambda (1+t)}{t}\bigg)^r = \frac{\log_\lambda(1+t)}{t (1+t)^k}\bigg(\frac{\log_\lambda(1+t)}{t}\bigg)^{r-1}(1+t)^k\\
&=\sum_{i=1}^\infty (-1)^{i+1} H_{i,\lambda}^{(k)} t^{i-1}\sum_{j=0}^\infty D_{j,\lambda}^{(r-1)} \frac{t^j}{j!}\sum_{l=0}^\infty \binom{k}{l}t^l\\
&=\sum_{i=0}^\infty (-1)^{i} H_{i+1,\lambda}^{(k)} t^{i}\sum_{m=0}^\infty \bigg(\sum_{j=0}^m \binom{m}{j} D_{j,\lambda}^{(r-1)}(k)_{m-j}\bigg)\frac{t^m}{m!} \\
&=\sum_{n=0}^\infty \bigg(\sum_{i=0}^n \sum_{j=0}^{n-i} (-1)^i \binom{n-i}{j}\frac{(k)_{n-i-j}}{(n-i)!} D_{j,\lambda}^{(r-1)} H_{i+1,\lambda}^{(k)}\bigg)t^n.
\end{split}
\end{equation}

Therefore, by comparing the coefficients on both sides of \eqref{eq25}, we obtain the following theorem.

\begin{theorem}
For $n, k \geq0$ and $r \geq1$, we have
\begin{equation*}
\begin{split}
D_{n,\lambda}^{(r)}=n!\sum_{i=0}^n \sum_{j=0}^{n-i}(-1)^i \binom{n-i}{j} \frac{(k)_{n-i-j}}{(n-i)!}D_{j,\lambda}^{(r-1)}H_{i+1,\lambda}^{(k)}.
\end{split}
\end{equation*}
\end{theorem}

%10
By \eqref{eq11}, we get
\begin{equation}\label{eq26}
\begin{split}
\sum_{n=1}^\infty H_{n,\lambda}t^n&=-\frac{\log_\lambda (1-t)}{1-t} = \frac{\log_\lambda(1-t)}{-t} \frac{t}{1-t} \\
&=\sum_{l=0}^\infty (-1)^l D_{l,\lambda} \frac{t^l}{l!}\sum_{j=1}^\infty t^{j} \\
&=\sum_{n=1}^\infty \bigg(\sum_{l=0}^{n-1}(-1)^l \frac{D_{l,\lambda}}{l!}\bigg)t^n.
\end{split}
\end{equation}
Thus, by comparing the coefficients on both sides of \eqref{eq26}, we get
\begin{equation}\label{eq27}
\begin{split}
H_{n,\lambda}=\sum_{l=0}^{n-1}(-1)^l \frac{D_{l,\lambda}}{l!}, \quad{(n\geq1)}.
\end{split}
\end{equation} \par
From \eqref{eq12}, we can derive the following.
\begin{equation}\label{eq28}
\begin{split}
\sum_{n=1}^\infty H_{n,\lambda}^{(r)}t^n&=-\frac{\log_\lambda (1-t)}{t} \frac{t}{(1-t)^r} \\
&=\sum_{l=0}^\infty D_{l,\lambda}(-1)^l\frac{t^l}{l!}\sum_{m=1}^\infty \binom{r+m-2}{m-1}t^m \\
&=\sum_{n=1}^\infty \bigg(\sum_{m=1}^n \binom{r+m-2}{r-1} \frac{D_{n-m,\lambda}}{(n-m)!}(-1)^{n-m}\bigg)t^n.
\end{split}
\end{equation}
Therefore, by \eqref{eq27} and \eqref{eq28}, we obtain the following theorem.

%11
\begin{theorem}
For $n \in \mathbb{N}$, we have
\begin{equation*}
\begin{split}
H_{n,\lambda}=\sum_{l=0}^{n-1}(-1)^l\frac{D_{l,\lambda}}{l!},\quad{(n\geq1)},
\end{split}
\end{equation*}

and

\begin{equation*}
\begin{split}
H_{n,\lambda}^{(r)}=\sum_{m=1}^n \binom{r+m-2}{r-1} \frac{D_{n-m,\lambda}}{(n-m)!}(-1)^{n-m}.
\end{split}
\end{equation*}
\end{theorem}

The degenerate derangements are defined by
\begin{equation}\label{eq29}
\begin{split}
\frac{1}{1-t}e_\lambda(-t)=\sum_{n=0}^\infty d_{n,\lambda}\frac{t^n}{n!}.
\end{split}
\end{equation}

Thus, we note that
\begin{equation*}
\begin{split}
d_{n,\lambda}=n!\sum_{k=0}^n(1)_{k,\lambda}\frac{(-1)^k}{k!}, \quad(n\geq0).
\end{split}
\end{equation*}

Now, we observe that
\begin{equation}\label{eq30}
\begin{split}
-\frac{\log_\lambda(1-t)}{(1-t)^r}e_\lambda(-t)&=\sum_{l=1}^\infty H_{l,\lambda}^{(r)}t^l\sum_{k=0}^\infty\frac{(1)_{k,\lambda}}{k!}(-1)^k t^k \\
&=\sum_{n=1}^\infty \bigg(\sum_{l=1}^nH_{l,\lambda}^{(r)}\frac{(1)_{n-l,\lambda}}{(n-l)!}(-1)^{n-l}\bigg)t^n.
\end{split}
\end{equation}

On the other hand, by \eqref{eq29}, we get
\begin{equation}\label{eq31}
\begin{split}
\frac{-\log_\lambda(1-t)}{(1-t)^r}e_\lambda(-t)&=-\frac{\log_\lambda(1-t)}{(1-t)^{r-1}}\frac{1}{1-t}e_\lambda(-t) \\
&=\sum_{l=1}^\infty H_{l,\lambda}^{(r-1)}t^l\sum_{k=0}^\infty d_{k,\lambda}\frac{t^k}{k!}=\sum_{n=1}^\infty \bigg(\sum_{l=1}^n H_{l,\lambda}^{(r-1)} \frac{d_{n-l,\lambda}}{(n-l)!}\bigg)t^n.
\end{split}
\end{equation}

%12
Therefore, by \eqref{eq30} and \eqref{eq31}, we obtain the following theorem.

\begin{theorem}

For $n \in \mathbb{N}$, we have
\begin{equation*}
\begin{split}
\sum_{l=1}^n H_{l,\lambda}^{(r)}\frac{(1)_{n-l,\lambda}}{(n-l)!}(-1)^{n-l}=\sum_{l=1}^n H_{l,\lambda}^{(r-1)} \frac{d_{n-l,\lambda}}{(n-l)!}.
\end{split}
\end{equation*}
\end{theorem}

We let $Y=\log_\lambda(1+t)$.
Then, for $ N \ge 1$, we have
\begin{equation}\label{eq32}
\begin{split}
\bigg(\frac{d}{dt}\bigg)^NY&=(\lambda-1)(\lambda-2)\cdots(\lambda-N+1)(1+t)^{\lambda-N} \\
&=\frac{N!}{\lambda}\binom{\lambda}{N}e_{\lambda}^{\lambda-N}(\log_\lambda(1+t))\\
&=\frac{N!}{\lambda} \binom{\lambda}{N}\sum_{k=0}^\infty (\lambda-N)_{k,\lambda}\frac{1}{k!}(\log_\lambda(1+t))^k \\
&=\frac{N!}{\lambda}\binom{\lambda}{N}\sum_{k=0}^\infty (\lambda-N)_{k,\lambda} \sum_{n=k}^\infty S_{1,\lambda}(n,k) \frac{t^n}{n!} \\
&=\sum_{n=0}^\infty \bigg(\frac{N!}{\lambda}\binom{\lambda}{N}\sum_{k=0}^nS_{1,\lambda}(n,k)(\lambda-N)_{k,\lambda}\bigg)\frac{t^n}{n!},
\end{split}
\end{equation}
where $N$ is a positive integer.\par
On the other hand, by \eqref{eq10}, we get
\begin{equation}\label{eq33}
\begin{split}
Y=\log_\lambda(1+t)=\frac{\log_\lambda(1+t)}{t}t=\sum_{n=1}^\infty nD_{n-1,\lambda}\frac{t^n}{n!}.
\end{split}
\end{equation}
Thus, by \eqref{eq33}, we get
\begin{equation}\label{eq34}
\begin{split}
\bigg(\frac{d}{dt}\bigg)^NY&=\sum_{n=N}^\infty n D_{n-1,\lambda}n(n-1)\cdots (n-N+1)\frac{t^{n-N}}{n!}\\
&=\sum_{n=0}^\infty (n+N)D_{n+N-1,\lambda}\frac{t^n}{n!}.
\end{split}
\end{equation}

Therefore, by \eqref{eq32} and \eqref{eq34}, we obtain the following theorem.

\begin{theorem}

For $N\in\mathbb{N}$ and $n\geq N-1$, we have
\begin{equation*}
\begin{split}
D_{n,\lambda}=\frac{N!}{n+1}\frac{1}{\lambda}\binom{\lambda}{N}\sum_{k=0}^{n-N+1}S_{1,\lambda}(n-N+1,k)(\lambda-N)_{k,\lambda}.
\end{split}
\end{equation*}
\end{theorem}

Next, we let $F=-\log_\lambda(1-t)$. Then, for $ N \ge 1$, we have
\begin{equation}\label{eq35}
\begin{split}
\bigg(\frac{d}{dt}\bigg)^NF&=(-1)^{N+1}(\lambda-1)(\lambda-2)\cdots(\lambda-N+1)(1-t)^{\lambda-N}\\
&=(-1)^{N+1}\frac{N!}{\lambda}\binom{\lambda}{N}e_\lambda^{\lambda-N}(\log_\lambda(1-t))\\
&=(-1)^{N+1}N!\frac{1}{\lambda}\binom{\lambda}{N}\sum_{k=0}^\infty (\lambda-N)_{k,\lambda}\frac{1}{k!}(\log_\lambda(1-t))^k\\
&=(-1)^{N+1}N!\frac{1}{\lambda}\binom{\lambda}{N}\sum_{k=0}^\infty (\lambda-N)_{k,\lambda}\sum_{n=k}^\infty S_{1,\lambda}(n,k)(-1)^n\frac{t^n}{n!}\\
&=\sum_{n=0}^\infty \bigg(N!\frac{1}{\lambda}\binom{\lambda}{N}\sum_{k=0}^n(-1)^{n-N-1}(\lambda-N)_{k,\lambda}S_{1,\lambda}(n,k)\bigg)\frac{t^n}{n!}.
\end{split}
\end{equation}

On the other hand, by \eqref{eq11}, we get
\begin{equation}\label{eq36}
F=-\log_\lambda(1-t)=-\frac{\log_\lambda(1-t)}{1-t}(1-t)
=\sum_{n=1}^\infty (H_{n,\lambda}-H_{n-1,\lambda})t^n.
\end{equation}

Thus, by \eqref{eq36} and for $N \ge 1$, we have
\begin{equation}\label{eq37}
\begin{split}
\bigg(\frac{d}{dt}\bigg)^NF&=\sum_{n=N}^\infty n(n-1)\cdots(n-N+1)(H_{n,\lambda}-H_{n-1,\lambda})t^{n-N}\\
&=\sum_{n=0}^\infty (n+N)(n+N-1)\cdots(n+1)(H_{n+N,\lambda}-H_{n+N-1,\lambda})t^n\\
&=\sum_{n=0}^\infty N!\binom{n+N}{N}(H_{n+N,\lambda}-H_{n+N-1,\lambda})t^n.
\end{split}
\end{equation}

Therefore, by \eqref{eq35} and \eqref{eq37}, we obtain the following theorem.
\begin{theorem}
For $N\in \mathbb{N}$ and $n\geq 0$, we have
\begin{equation*}
\begin{split}
\frac{1}{n!}\frac{1}{\lambda} \binom{\lambda}{N} \sum_{k=0}^n (-1)^{n-N-1}(\lambda-N)_{k,\lambda}S_{1,\lambda}(n,k)=\binom{n+N}{N}(H_{n+N,\lambda}-H_{n+N-1,\lambda}).
\end{split}
\end{equation*}

\end{theorem}

By Theorem 9 and \eqref{eq06}, we get
\begin{equation}\label{eq38}
\begin{split}
\frac{1}{n!}&\sum_{k=0}^n(-1)^{n-N-1}(\lambda-N)_{k,\lambda}S_{1,\lambda}(n,k)=\frac{\binom{n+N}{N}}{\frac{1}{\lambda}\binom{\lambda}{N}}(H_{n+N,\lambda}-H_{n+N-1,\lambda})\\
&=\frac{\binom{n+N}{N}}{\frac{1}{\lambda}\binom{\lambda}{N}}\frac{1}{\lambda}\binom{\lambda}{n+N}(-1)^{n+N-1}=(-1)^{n+N-1}\frac{\binom{\lambda}{N+n}}{\binom{\lambda}{N}}\binom{n+N}{N}.
\end{split}
\end{equation}

Therefore, by \eqref{eq38}, we obtain the following corollary.
\begin{corollary}
For $n\geq0$ and $N\in\mathbb{N}$, we have
\begin{equation*}
\begin{split}
\frac{1}{n!}\sum_{k=0}^n(\lambda-N)_{k,\lambda}S_{1,\lambda}(n,k)=\frac{\binom{\lambda}{n+N}}{\binom{\lambda}{N}}\binom{n+N}{N}.
\end{split}
\end{equation*}
\end{corollary}

%16
\begin{remark}
From Corollary 10 and letting $\lambda \rightarrow 0$, we obtain
\begin{equation*}
(-1)^{n}\frac{N}{n+N}\binom{n+N}{N}=\frac{1}{n!}\sum_{k=0}^n(-1)^{k}N^kS_1(n,k).
\end{equation*}
\end{remark}

\begin{remark}
Recently, on the Daehee numbers and their related topics various studies have been conducted by several researchers. Interested readers may refer to \cite{1, 12, 13, 14, 15, 17, 18}.
\end{remark}

\hspace{5cm}

\section{conclusion}
Many different tools have been used in the explorations for degenerate versions of some special numbers and polynomials, which include generating functions, combinatorial methods, umbral calculus, $p$-adic analysis, differential equations, probability theory, operator theory, special functions and analytic number theory (\text{see \cite{5, 6, 7, 8, 9, 10, 11, 12}} and the references therein). In this paper, we used the elementary methods of generating functions in order to study the degenerate harmonic and degenerate hyperharmonic numbers. Some properties, recurrence relations and identities relating to those numbers were derived in connection with the degenerate Stirling numbers of the first kind, the degenerate Daehee numbers and the degenerate derangement. \par
We would like to continue to investigate various degenerate versions of certain special numbers and polynomials, especially their applications to physics, science and engineering.

\vspace{0.2in}

\noindent{\bf{Acknowledgments}} \\
%The author would like to thank the referees for the detailed and valuable comments that helped
%improve the original manuscript in its present form.% Also,
%%
The authors thank Jangjeon Institute for Mathematical Sciences for the support of this research.

\vspace{0.1in}

\noindent{\bf {Availability of data and material}} \\
Not applicable.

\vspace{0.1in}

\noindent{\bf{Funding}} \\
This work was supported by the Basic Science Research Program, the National
               Research Foundation of Korea,
               (NRF-2021R1F1A1050151).
\vspace{0.1in}

\noindent{\bf{Ethics approval and consent to participate}} \\
All authors declare that there is no ethical problem in the production of this paper.

\vspace{0.1in}

\noindent{\bf {Competing interests}} \\
All authors declare no conflict of interest.

\vspace{0.1in}

\noindent{\bf{Consent for publication}} \\
All authors want to publish this paper in this journal.

\vspace{0.1in}

\noindent{\bf{Author' Contributions}}\\
All authors read and approved the final manuscript.

\end{document}